\documentclass[12pt]{article}
\usepackage{latexsym}
\usepackage[dvips]{graphicx}
\usepackage{amssymb, amsmath}
\usepackage{bm}
\usepackage{fleqn}
\usepackage{subfigure}

\renewcommand{\labelenumi}{(\arabic{enumi})}

\def\disp{\displaystyle } 
\def\G{\Gamma}
\def\P{{\cal P}}
\def\S{{\cal S}}
\def\R{\mathbb{R}}
\def\Z{\mathbb{Z}}

\def\veca{\bm{a}}
\def\vecb{\bm{b}}
\def\vecp{\bm{p}}
\def\vecx{\bm{x}}

\def\vecT{\bm{T}}
\def\vecN{\bm{N}}
\def\sgn{\mathop\mathrm{sgn}}
\def\Image{\mathop\mathrm{Image}}

\def\<{\langle}
\def\>{\rangle}

\def\eps{\varepsilon}
\begin{document}
\title{\vspace{-2truecm}
Application of a curvature adjusted method in image segmentation\footnote{
MB, MK, PP and SY are supported by Czech Technical University in Prague, 
Faculty of Nuclear Sciences and Physical Engineering, within the 
Jind\v rich Ne\v cas Center for 
Mathematical Modeling (Project of the Czech Ministry of Education, 
Youth and Sports LC 06052).
TT is partially supported by Grant-in-Aid for Scientific Research (No.17540125) 
by Japan Society for the Promotion of Science. 
DS was supported by APVV-0247-06. 
}}
\author{\normalsize 
$\mbox{M. Bene\v{s}}^{\mathrm{a})}$, 
$\mbox{M. Kimura}^{\mathrm{a}), \mathrm{b})}$, 
$\mbox{P. Pau\v{s}}^{\mathrm{a})}$, 
$\mbox{D. \v{S}ev\v{c}ovi\v{c}}^{\mathrm{c})}$, 
$\mbox{T. Tsujikawa}^{\mathrm{d})}$ 
\ and \ 
$\mbox{S. Yazaki}^{\mathrm{a}), \mathrm{d})}$}
\date{}
\maketitle
\vspace{-1truecm}
\begin{center}
\begin{minipage}[h]{0.95\linewidth}{\small
a) Faculty of Nuclear Sciences and Physical Engineering, 
Czech Technical University in Prague, 
Trojanova 13, 120 00 Prague, Czech Republic. 
\\\small
b) Faculty of Mathematics, 
Kyushu University, 
6-10-1 Hakozaki, 
Fukuoka 812-8581, Japan. 
\\\small
c) Institute of Applied Mathematics, 
Faculty of Mathematics, Physics and Informatics, 
Comenius University, 
842 48 Bratislava, Slovak Republic. 
\\\small
d) Faculty of Engineering, 
University of Miyazaki, 
1-1 Gakuen Kibanadai Nishi, 
Miyazaki 889-2192, Japan. 
}\end{minipage}
\end{center}
\begin{abstract}
This article deals with 
flow of plane curves driven by the curvature and external force.
We make use of such a geometric flow for the purpose of image segmentation. 
A parametric model for evolving curves with uniform and curvature adjusted redistribution 
of grid points will be described and compared. 
\end{abstract}

\section{\normalsize Introduction}

Finding shapes in images and their transformation to some usable form are very 
important and attractive tasks in numerical computation. 
For such an image segmentation, various numerical methods have been proposed. 
In this article, we will discuss 2D image segmentation 
by using moving plane curves which evolve according to the law: 
\begin{equation}\label{v=k+F}
\mbox{normal velocity}=\mbox{curvature}+\mbox{external force}. 
\end{equation}
Under suitable setting of an external force and appropriate choice of 
numerical parameters, 
segmentation of an edge of the given object image will be obtained as a 
stable stationary solution. 
As it will be obvious from numerical examples in the following sections,
one can obtain a good approximation of the image edge after 
sufficiently long evolution of planer curves.
The CPU time is very short, varying from 1 to 10 seconds order. 

Evolving curves can be treated numerically as well as mathematically in several ways: 
for instance, 
direct approach based methods (Lagrangian method), 
level-set methods, phase-field methods, etc. 
Our method is based on a direct approach. 
In the direct approach, however, 
evolving curves are approximated by polygons, the vertices are tracked with discrete time. 
This type of tracking methods often may produce unstable computational results. 
For instance, 
vertices may accumulate somewhere and may be very sparse elsewhere. 
Hence, a more stable numerical tracking scheme for evolving curves is required. 
To avoid these undesired phenomena, 
some kind of points redistribution method 
have been extensively studied by many authors. 
One of the useful ways is redistribution of points uniformly along the curve by means of the tangential velocity. 
In the series of papers \cite{MikulaS2001, MikulaS2004a, MikulaS2004b, MikulaS2006},
K. Mikula and D. \v Sev\v covi\v c 
proposed asymptotically uniform redistribution method for a general 
geometric equation (\ref{v=k+F}), and,
in particular, they applied it to image segmentation. 
We also refer to \cite{SrikrishnanCRS2007} for another approach to image segmentation
of moving objects,
where asymptotically uniform redistribution method 
has been applied for a geometric equation (\ref{v=k+F}) with
the external force including 
the probability that the curve belongs to the background or 
target of the segmented sequence of images. 
For image segmentation of indirect approach, 
we refer, e.g., \cite{BenesCM2004, FrolkovicMPS20XX} for the readers. 

Besides the uniform redistribution method, 
a curvature adjusted method was proposed by D. \v Sev\v covi\v c and S. Yazaki 
in \cite{SevcovicY20XXa}, 
in which one can find a short history of redistribution methods. 
Their method provides adequate redistribution of grid points 
depending on the absolute value of curvature. 
The idea comes from the tangential 
redistribution in the case of a crystalline curvature flow \cite{Yazaki20XXa}. 
The advantage of this method is that few points are needed to get sharp corners on edges
in comparison with the uniform redistribution method.

Organization of the present paper is as follows. 
In the next section, the evolution equations will be introduced. 
In Section 3, a numerical scheme is proposed and analyzed.
In the last section, computational experiments are presented. 

\section{\normalsize Evolution equations of curves}

Our method is based on the direct parametric 
(Lagrangian) approach. 
According to \cite{SevcovicY20XXa}, 
we introduce the evolution equation as follows. 
We consider an embedded and closed plane curve $\G$ which is 
parametrized counterclockwise 
by a smooth periodic function $\vecx:\ \R/\Z\supset[0, 1]\to\R^2$ such that 
$\G=\Image(\vecx)=\{\vecx(u);\ u\in[0, 1]\}$ and
$|\partial_u\vecx|>0$. 
Here and after, 
the circle $\R/\Z$ is represented by $[0, 1]$ with periodic condition 
$\vecx(0)=\vecx(1)$, 
and we define $\partial_\xi{\sf F}=\partial{\sf F}/\partial\xi$, 
and $|\veca|=\sqrt{\veca.\veca}$ where 
$\veca.\vecb$ is the Euclidean inner product between vectors $\veca$ and $\vecb$. 
The unit tangent vector can be defined as 
$\vecT=\partial_u\vecx/|\partial_u\vecx|=\partial_s\vecx$, 
where $s$ is the arc-length parameter and $ds=|\partial_u\vecx|du$, and 
the unit inward normal vector is defined in such a way that 
$\vecT\wedge\vecN=1$, where $\veca\wedge\vecb=\det(\veca, \vecb)$ for 
two dimensional column vectors $\veca$, $\vecb$. 
The signed curvature in the direction $\vecN$ is denoted by $k$,
e.g.,  the sign of $k$ is plus for strictly convex curves.
Let $\nu$ be the angle of $\vecT$, i.e., $\vecT=(\cos\nu, \sin\nu)$ and 
$\vecN=(-\sin\nu, \cos\nu)$. 
The geometric evolution law (\ref{v=k+F}) of plane curves is stated as follows: 
For a given initial curve $\G^0=\Image(\vecx^0)=\G$, find a family of curves 
$\{\G^t\}_{t\geq 0}$, 
$\G^t=\{\vecx(u, t);\ u\in[0, 1]\}$ starting from
$\vecx(u, 0)=\vecx^0(u)$ for $u\in[0, 1]$ and evolving according to the 
geometric equation
\begin{equation}\label{eq:v=beta}
v = k + F(\vecx),
\end{equation}
where $v$ is the normal velocity in the inward normal direction, and $F$ is a 
given external force defined everywhere on the plane. 
The evolution equation of a position vector $\vecx$ is given as: 
\begin{equation}\label{eq:basic_eq}
\partial_t\vecx=(k+F)\vecN+\alpha\vecT, \quad
\vecx(\cdot, 0)=\vecx^0(\cdot). 
\end{equation}
Here $\alpha$ is the tangential component of the velocity vector. 
Notice that the presence of an arbitrary
tangential component does not affect the shape of curve
(see \cite[Proposition 2.4]{EpsteinG1987}). 
Also, $\partial_t{\sf F}$ means $\partial_t{\sf F}(u, t)$, and 
$\partial_s{\sf F}$ is given in terms of $u$: 
$\partial_s{\sf F}=|\partial_u\vecx|^{-1}\partial_u{\sf F}$. 

The external force $F(\vecx)$ plays a 
very important role in the application to the image segmentation, 
because it describes the image to be segmented.
In our approach, its value is given by the intensity of the image on each pixel. 
The formula how to calculate the value of $F$ from the picture will 
be mentioned in (\ref{eq:color}). 

According to \cite{MikulaS2001}, 
equation (\ref{eq:basic_eq}) can be written as the following PDE:
\begin{equation}
\partial_t\vecx=\partial_s^2\vecx+F\vecN+\alpha\partial_s\vecx \,.
\label{eq:equation_position}
\end{equation}
A solution $\vecx$ is subject to the initial condition $\vecx(\cdot, 0)=\vecx^0(\cdot)$. 
The curve is closed, so we have to impose periodic boundary 
conditions of $\vecx$ for $u\in[0,1]$ and $\nu(1,t)=\nu(0,t)+2\pi$.
The term $\alpha$ is responsible for the tangential redistribution.
According to \cite{SevcovicY20XXa, SevcovicY20XXb},
the curvature adjusted tangential velocity has
been proposed in the form:
\begin{equation}\label{eq:tangential_velocity}
\partial_s(\varphi(k)\alpha)
=f-\frac{\varphi(k)}{\<\varphi(k)\>}\<f\>
+\omega\left(\frac{L^t}{|\partial_u\vecx|}\<\varphi(k)\>-\varphi(k)\right), 
\end{equation}
where $L^t$ is the curve length in time $t$, 
$\omega$ is a given positive constant, and 
\begin{equation}\label{eq:phi}
\varphi(k) = 1 - \varepsilon + \varepsilon \sqrt{1-\varepsilon +\varepsilon k^2},
\end{equation}
\[
f=\varphi(k)k (k+F)-\varphi'(k)\left(\partial_s^2 k+ \partial_s^2 F +k^2(k+F) \right), \quad
\varphi'(k)=\frac{d}{dk}\varphi(k), 
\]
\[
\<{\sf F}(\cdot, t)\>=\frac{1}{L^t}\int_{\G^t}{\sf F}(s, t)\,ds. 
\]
Equation (\ref{eq:tangential_velocity}) is designed for the exponential convergence of an 
extended relative local length 
$\lim_{t\to+\infty}|\partial_u\vecx|\varphi(k)/(L^t\<\varphi(k)\>)=1$, 
with the exponent $-\omega t$. 
The function $\varphi(k)$ given by (\ref{eq:phi}) is very important
because it controls the redistribution of grid points. 
Note that $\varphi(k)\to 1$ if $\eps\to 0^{+}$, and 
$\varphi(k)\to|k|$ if $\eps\to 1^{-}$. 
The function $\varphi(k)\equiv 1$ produces the uniform redistribution of points for $\omega=0$ 
and asymptotically uniform redistribution of points for $\omega>0$ 
(see \cite{MikulaS2001, MikulaS2004a, MikulaS2004b, MikulaS2006}), 
and the function $\varphi(k)=|k|$ is used implicitly for the crystalline curvature flow \cite{Yazaki20XXa}. 
By choosing $\eps\in(0, 1)$,
we obtain curvature adjusted redistribution \cite{SevcovicY20XXa}. 

To get unique solution $\alpha$, 
we assume the following additional condition for $\alpha$:
\begin{equation}\label{eq:tangential_velocity_av=0}
\<\alpha(\cdot, t)\>=0. 
\end{equation}

The details can be found in the forthcoming paper \cite{SevcovicY20XXb}. 

\section{\normalsize Numerical scheme}

We follow the numerical scheme developed in \cite{SevcovicY20XXa, SevcovicY20XXb}.
For a given initial $N$-sided polygonal curve $\P^0=\bigcup_{i=1}^N\S_i^0$, 
let us find a family of $N$-sided polygonal curves $\{\P^j\}_{j=1, 2, \ldots}$, 
$\P^j=\bigcup_{i=1}^N\S_i^j$, where 
$\S_i^j=[\vecx_{i-1}^j, \vecx_{i}^j]$ is the $i$th edge with 
$\vecx_{0}^j=\vecx_N^j$ for $j=0, 1, 2, \ldots$. 
The initial polygon $\P^0$ is an approximation of $\G^0$.
We construct $\P^j$ as an approximation of $\G^t$ at the $j$th time $t=t_j$, 
where $t_0=0$ and $t_j=\sum_{l=0}^{j-1}\tau_l$ ($j=1, 2, \ldots$) with 
the adaptive time increments $\tau_l>0$ for $l=0,1,\ldots$. 
The updated curve $\{\P^{j+1}\}$ is determined from the data $\{\P^j\}$ at the previous time step by 
using discretization in space and time of a closed system of PDEs 
(\ref{eq:equation_position}) and (\ref{eq:tangential_velocity}). 

The way of discretization is similar to \cite{MikulaS2004b}, it means that
our scheme is semi-implicit scheme and the discretization is based on 
the flowing finite volume approach from
\cite{MikulaS2001}. See also \cite{MinarikKMB2004, MinarikKM2005}. 
Discrete quantities 
$\alpha_i^j$, $k_i^j$, $\nu_i^j$, $\vecx_i^j$, $r_i^j=|\S_i^j|$ of $\P^j$ 
($i=1, 2, \ldots, N$) are splitted into two categories: 
$k_i^j$, $\nu_i^j$ and $r_i^j$ take constant value on the 
$i$th edge $\S_i^j$ (the finite volume),
whereas $\alpha_i^j$ and $\vecx_i^j$ are defined at the $i$th vertex $\vecx_i^j$ and take 
constant values on the $i$th dual edge 
$\S_i^{*j}=[\vecx_i^{*j}, \vecx_{i+1}^{*j}]$ (the corresponding dual volume)
where $\vecx_i^{*j}=(\vecx_i^j+\vecx_{i-1}^j)/2$. 
The duality is defined in such a way that 
${\sf F}_i^{*j}=({\sf F}_i^{j}+{\sf F}_{i+1}^{j})/2$ 
(where ${\sf F}=k, \nu, r$) 
take constant value on $\S_i^{*j}$, and 
$\alpha_i^{*j}=(\alpha_i^{j}+\alpha_{i-1}^{j})/2$ 
takes constant value on $\S_i^j$. 
Hereafter all quantities except $\{\nu_i^j\}$ satisfy the periodic boundary condition
${\sf F}_0={\sf F}_N$, ${\sf F}_{N+1}={\sf F}_1$. 
We will use the following notation: 
$(\partial_\tau{\sf F})_i^j=({\sf F}_i^{j+1}-{\sf F}_i^j)/\tau$, 
${\sf F}_{\rm min}=\min_{1\leq i\leq N}{\sf F}_i$, 
$|{\sf F}|_{\rm max}=\max_{1\leq i\leq N}|{\sf F}_i|$, 
$\<{\sf F}^j\>=\sum_{l=1}^N{\sf F}_l^jr_l^j/L^j$, 
$L^j=\sum_{l=1}^Nr_l^j$ is the total length of $\P^j$, 
$(\partial_s{\sf F})_i^j=({\sf F}_{i}^{j}-{\sf F}_{i-1}^{j})/r_{i}^{j}$, 
$(\partial_s{\sf F}^*)_i^j=({\sf F}_{i}^{*j}-{\sf F}_{i-1}^{*j})/r_{i}^{j}$, 
$(\partial_{s*}{\sf F})_i^j=({\sf F}_{i+1}^{j}-{\sf F}_{i}^{j})/r_{i}^{*j}$, 
$(\partial_{s*}{\sf F}^*)_i^j=({\sf F}_{i+1}^{*j}-{\sf F}_{*i}^{j})/r_{i}^{*j}$; and 
$(\forall i)$ means $(i=1, 2, \ldots, N)$, whereas
$(\forall' i)$ means $(i=2, 3, \ldots, N)$.

At first, the initial values for $\{r_i^0\}$,$\{k_i^0\}$,$\{\nu_i^0\}$ are computed
from the initial data according to (\ref{discrete_quant}).
\begin{equation}\label{discrete_quant}
\begin{array}{ll}\disp
(\mathrm{i}) &
r_i^j=|\vecp_i|, \quad \vecp_i=(p_{i_1}, p_{i_2})=\vecx_i^j-\vecx_{i-1}^j\quad 
(\forall i), \\[5pt]\disp
(\mathrm{ii}) & \disp
k_i^j=\frac{1}{2r_i^j}\sgn(\vecp_{i-1}\wedge\vecp_{i+1})
\arccos\left(\frac{\vecp_{i-1}.\vecp_{i+1}}{r_{i-1}^jr_{i+1}^j}\right)\quad 
(\forall i), \\[5pt]\disp
(\mathrm{iii}) &
\tilde{\nu}_i=\left\{\begin{array}{@{}ll}\disp
\arccos(p_{i_1}/r_i^j) & (p_{i_2}\geq 0)\\\disp
2\pi-\arccos(p_{i_1}/r_i^j) & (p_{i_2}<0)
\end{array}\right.\quad
(\forall i), \\[5pt]\disp
(\mathrm{iv}) &
\nu_1^j=\tilde{\nu}_1, \\[5pt]\disp
(\mathrm{v}) &
\nu_i^j=\left\{\begin{array}{@{}ll}\disp
\tilde{\nu}_i\pm 2\pi & (|\tilde{\nu}_i-\nu_{i-1}^j|>|\tilde{\nu}_i\pm 2\pi-\nu_{i-1}^j|)\\\disp
\tilde{\nu}_i & (\mbox{otherwise})\ 
\end{array}\right.\quad
(\forall' i), \\[5pt]\disp
(\mathrm{vi}) &
\nu_{N+1}^j=\left\{\begin{array}{@{}ll}\disp
\tilde{\nu}_1\pm 2\pi & (|\tilde{\nu}_1-\nu_N^j|>|\tilde{\nu}_1\pm 2\pi-\nu_N^j|), \\\disp
\tilde{\nu}_1 & (\mbox{otherwise}), 
\end{array}\right.\\[15pt]\disp
(\mathrm{vii}) &
\nu_0^j=\nu_1^j-(\nu_{N+1}^j-\nu_N^j). 
\end{array}
\end{equation}

Next we discretize equation (\ref{eq:tangential_velocity}) for redistribution term
$\alpha$ and compute its values.
By taking discrete time stepping, we obtain the following closed system of
finite difference equations of tangential velocity $\{\alpha_i^{j+1}\}$:
\begin{equation}\label{eq:FDEs-tangential_velocity}
\varphi(k_i^{*j})\alpha_i^{j+1}=
\varphi(k_{i-1}^{*j})\alpha_{i-1}^{j+1}
+\psi_i^j\quad (\forall' i), 
\end{equation}
where
\[
\psi_i^j=
f_i^jr_i^j
-\frac{\varphi(k_i^j)}{\<\varphi(k^j)\>}\<f^j\>r_i^j
+\omega\left(\frac{L^j}{N}\<\varphi(k^j)\>-\varphi(k_i^j)r_i^j\right),
\]
\[
f_i^j
=\varphi(k_i^j)k_i^j(k_i^j + F_i^{*j})
\]
\[
\quad -\frac{\varphi'(k_i^j)}{r_i^{j}}\left(\frac{k_{i+1}^j-k_i^j}{r^{*j}_i} - \frac{k_{i}^j-k_{i-1}^j}{r^{*j}_{i-1}}+
\frac{F_{i+1}^j-F_i^j}{r^{*j}_i} - \frac{F_{i}^j-F_{i-1}^j}{r^{*j}_{i-1}}\right),
\]
where 
$\omega>0$ is a given constant (in the following experiment, $\omega=50000$ will be used) and 
$F_i^{*j}$ = $F(\vecx_i^{*j})$. 
From (\ref{eq:FDEs-tangential_velocity}) we obtain
\begin{equation}\label{eq:alpha}
\varphi(k_i^{*j})\alpha_i^{j+1}
=\varphi(k_1^{*j})\alpha_1^{j+1}
+\Psi_i^j, \quad
\Psi_i^j=\sum_{l=2}^i\psi_l^j 
\quad (\forall' i). 
\end{equation}
On the other hand, the equality
$\sum_{i=1}^Nr_i^{*j}\alpha_i^{j+1}=0$ follows from $\<\alpha^{*j+1}\>=0$ which is 
a discretization of (\ref{eq:tangential_velocity_av=0}).
Therefore $\{\alpha_i^{j+1}\}_{i=1}^{N}$ can be determined uniquely from (\ref{eq:alpha}) and 
\[
\varphi(k_1^{*j})\alpha_1^{j+1}
=-\frac{\sum_{i=2}^Nr_i^{*j}\Psi_i^j/\varphi(k_i^{*j})}{\sum_{i=1}^Nr_i^{*j}/\varphi(k_i^{*j})}. 
\]

Finally, we discretize equation (\ref{eq:basic_eq}) 
for the position vector $\vecx$. According to \cite{SevcovicY20XXa, SevcovicY20XXb},
we obtain the following closed semi-implicit system:
\[
(\partial_\tau\vecx)_i^j
=(\partial_{s*}(\partial_s\vecx))_i^{j+1}
+\alpha_i^{j+1}(\partial_{s*}\vecx^*)_i^{j+1}
+F(\vecx_i^j)\vecN(\nu_i^{*j}). 
\] 
This system is equivalent to the following linear system for position vectors $\vecx_i^{j+1}$
subject to periodic boundary conditions: 
\begin{equation}\label{eq:FDEs-equation_position}
-a_i^{j+\frac12}\tau \vecx_{i-1}^{j+1}
+(1+b_i^{j+\frac12}\tau)\vecx_{i}^{j+1}
-c_i^{j+\frac12}\tau \vecx_{i+1}^{j+1}
=\vecx_{i}^{j}+F(\vecx_i^j)\vecN(\nu_i^{*j})\tau\quad (\forall i),
\end{equation}
where $b_i^{j+\frac12}=a_i^{j+\frac12}+c_i^{j+\frac12}$, and 
\[
a_i^{j+\frac12}=\frac{1}{r_{i}^{*j}}
\left(\frac{1}{r_{i}^{j}}
-\frac{\alpha_i^{j+1}}{2}\right), \quad
c_i^{j+\frac12}=\frac{1}{r_{i}^{*j}}
\left(\frac{1}{r_{i}^{j}}
+\frac{\alpha_i^{j+1}}{2}\right). 
\]
The above linear system is a diagonally dominant and solvable 
if we choose the following time step:
\[
\tau
=\frac{1}{\eta(1+\lambda)}, \quad
\eta=\frac{4}{r_{\rm min}^{j}}\left(
\frac{1}{r_{\rm min}^{j}}
+\frac{|\alpha^{j+1}|_{\rm max}}{2}\right), \quad
\lambda>0. 
\]

The transformation formula for mapping pixel colors to $F$ will be given in the next section. 

Our scheme is described in the following algorithm.\\
\begin{enumerate}
\renewcommand{\labelenumi}{\arabic{enumi}.}
\item Set $j:=0$. Start up with a sufficiently large initial closed curve $\{\P^0\}$
given by its vertices $\{\vecx_i^0\}_{i=1}^N$. 
\item Compute $\{r_i^j\}_{i=1}^N$, $\{k_i^j\}_{i=1}^N$, $\{\nu_i^j\}_{i=1}^N$ 
by means of equations (\ref{discrete_quant}).
\item Get the external force $\{F_i^{j}\}_{i=1}^N$ 
and $\{F_i^{*j}\}_{i=1}^N$ from the intensity function $I(\vecx)$.
\item Compute $\{\alpha_i^{j+1}\}_{i=1}^N$ by solving 
equation (\ref{eq:FDEs-tangential_velocity}).
\item Compute $\{\vecx_i^{j+1}\}_{i=1}^N$ by solving 
equation (\ref{eq:FDEs-equation_position}).
\item Set $j:=j+1$ and go to step 2. 
\end{enumerate}
In our numerical experiment, 
we always obtain almost stationary shape $\P^j$, if $j$ is large enough. 
An appropriate choice of stopping condition for the above algorithm 
is studied in our on-going research. 

\section{\normalsize Experimental comparison between the uniform and curvature adjusted redistribution methods}

\begin{figure}[h]
\begin{center}
\includegraphics[width=0.4\textwidth]{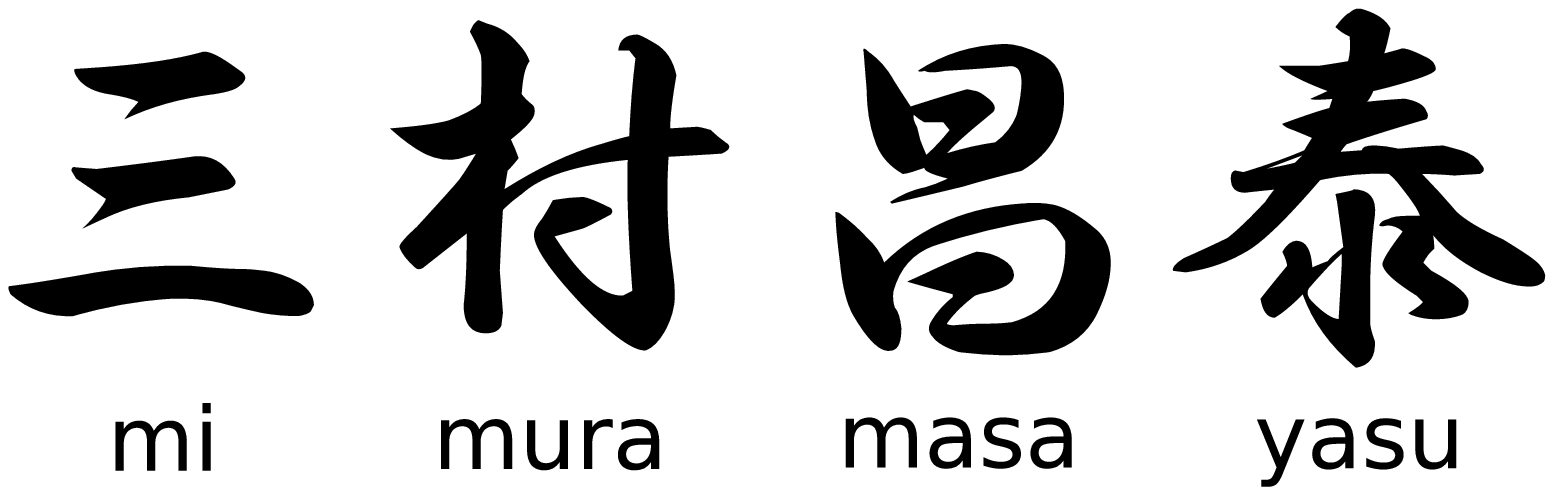}
\end{center}
\caption{The name of Professor Masayasu Mimura written in kanji. }
\label{fig:mimura}
\end{figure}

\begin{figure}[h]
\begin{center}
\begin{picture}(175,175)
\put(-70,19){\includegraphics[width=0.33\textwidth]{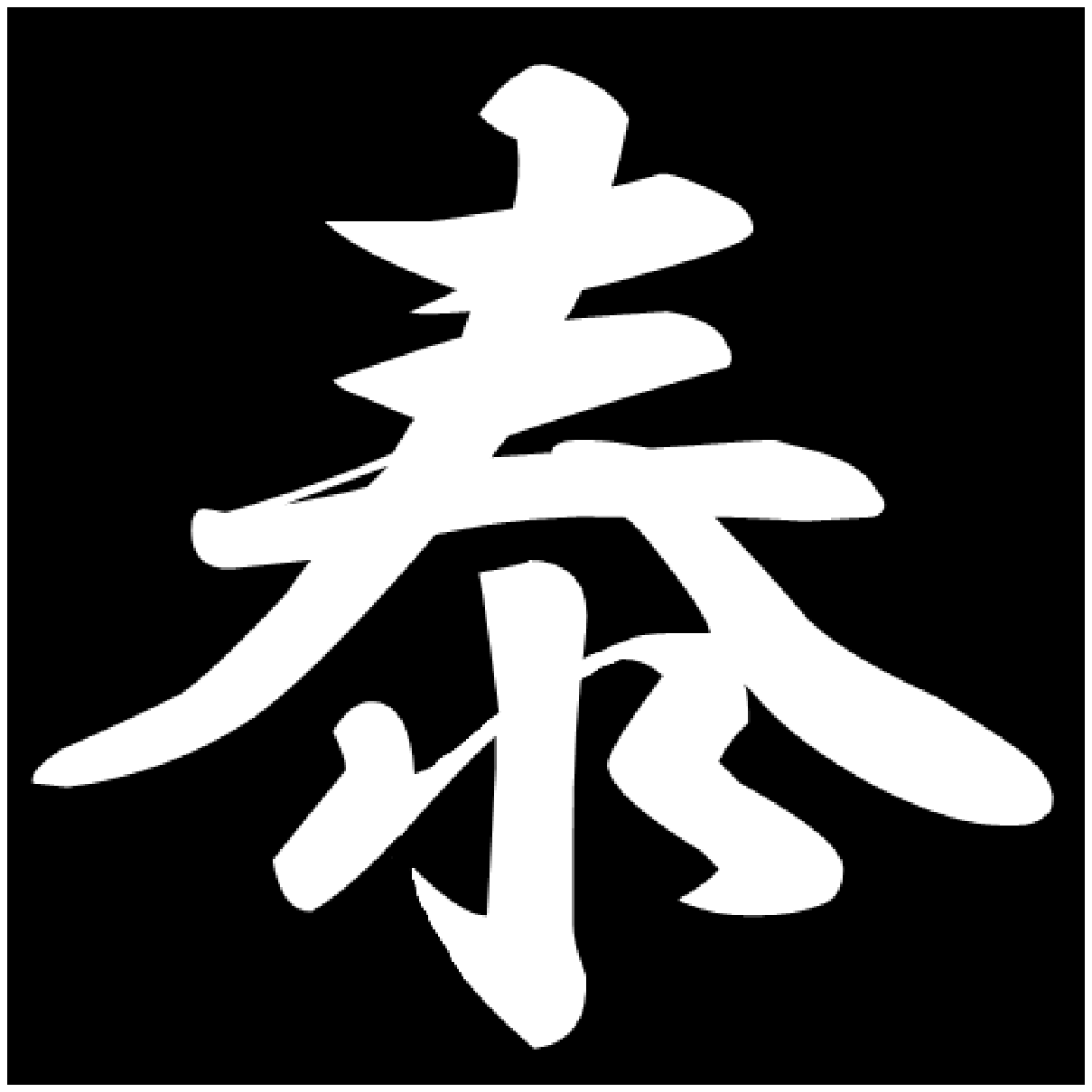}}
\put(100,0){\includegraphics[width=0.4\textwidth]{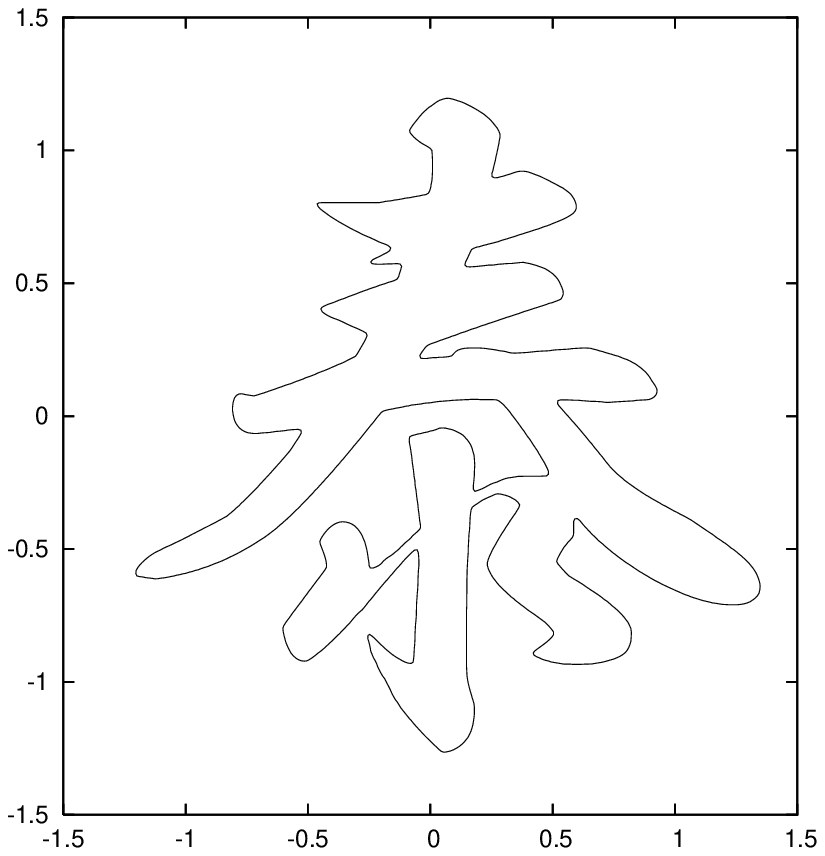}}
\end{picture}
\end{center}
\caption{The original bitmap image (left) 
and the very accurate image segmentation at $t=0.025$.
We chose $N=1500$, $\varepsilon = 0$, 
and $F\in (-100,100)$ (right).}
\label{fig:yasu_bmp}
\end{figure}

\begin{figure}[h]
\begin{center}
\begin{picture}(175,175)
\put(-70,20){
 \includegraphics[width=0.33\textwidth]{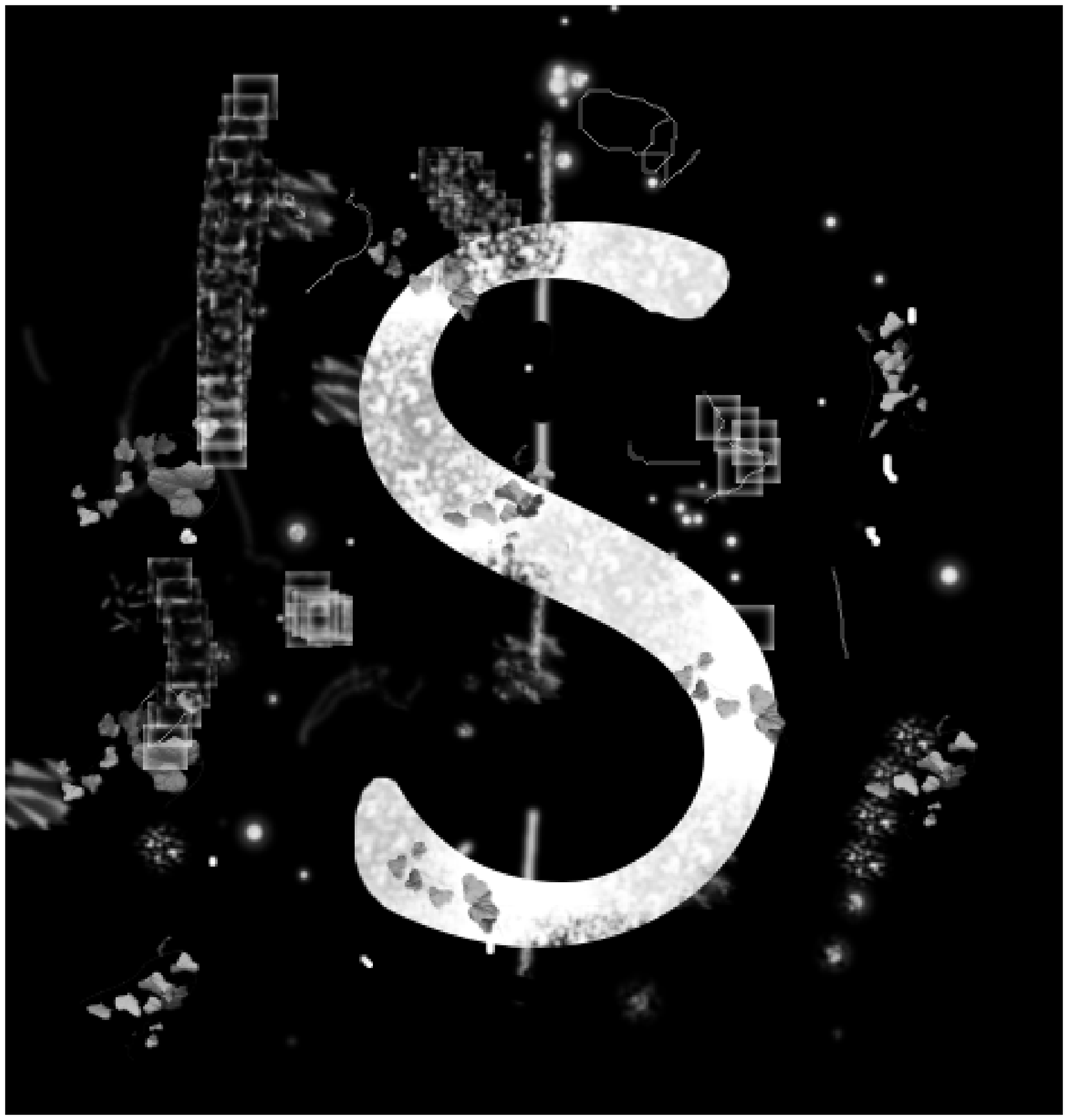}
 }
\put(100,0){
\includegraphics[width=0.4\textwidth]{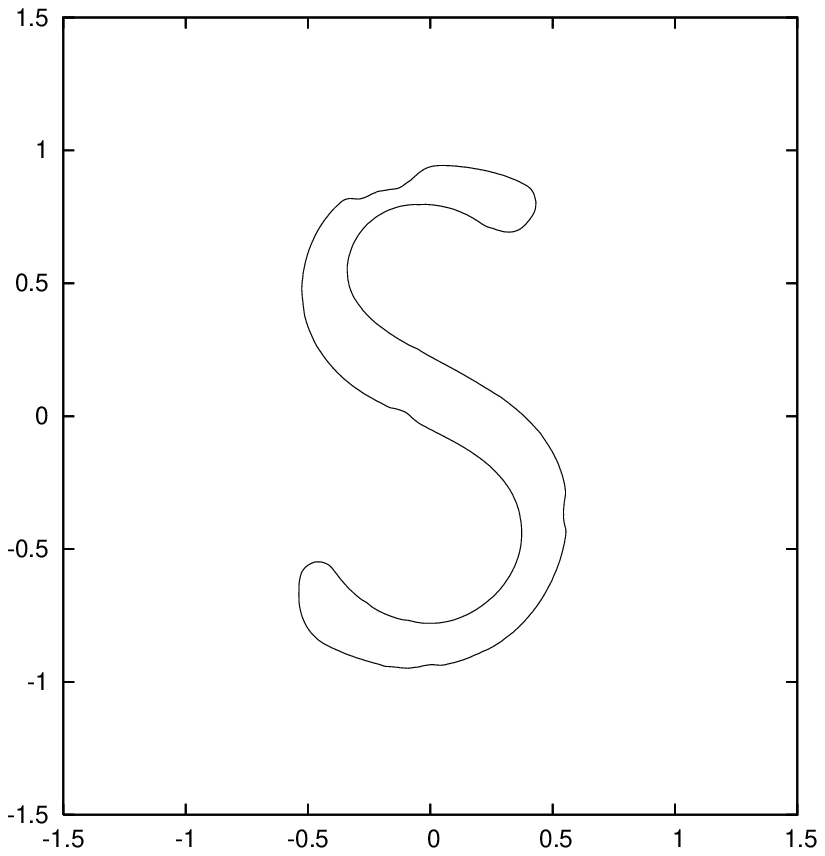}
}
\end{picture}
\end{center}
\caption{The original bitmap image (left) and a very accurate 
image segmentation at $t=0.1$. We chose $N=800$, $\varepsilon = 0$,
and $F\in (-30,35)$ (right).}
\label{fig:s_bmp}
\end{figure}

\begin{figure}[h]
 \begin{center}
 \begin{picture}(175,175)
 \put(-70,13){
 \includegraphics[width=0.36\textwidth]{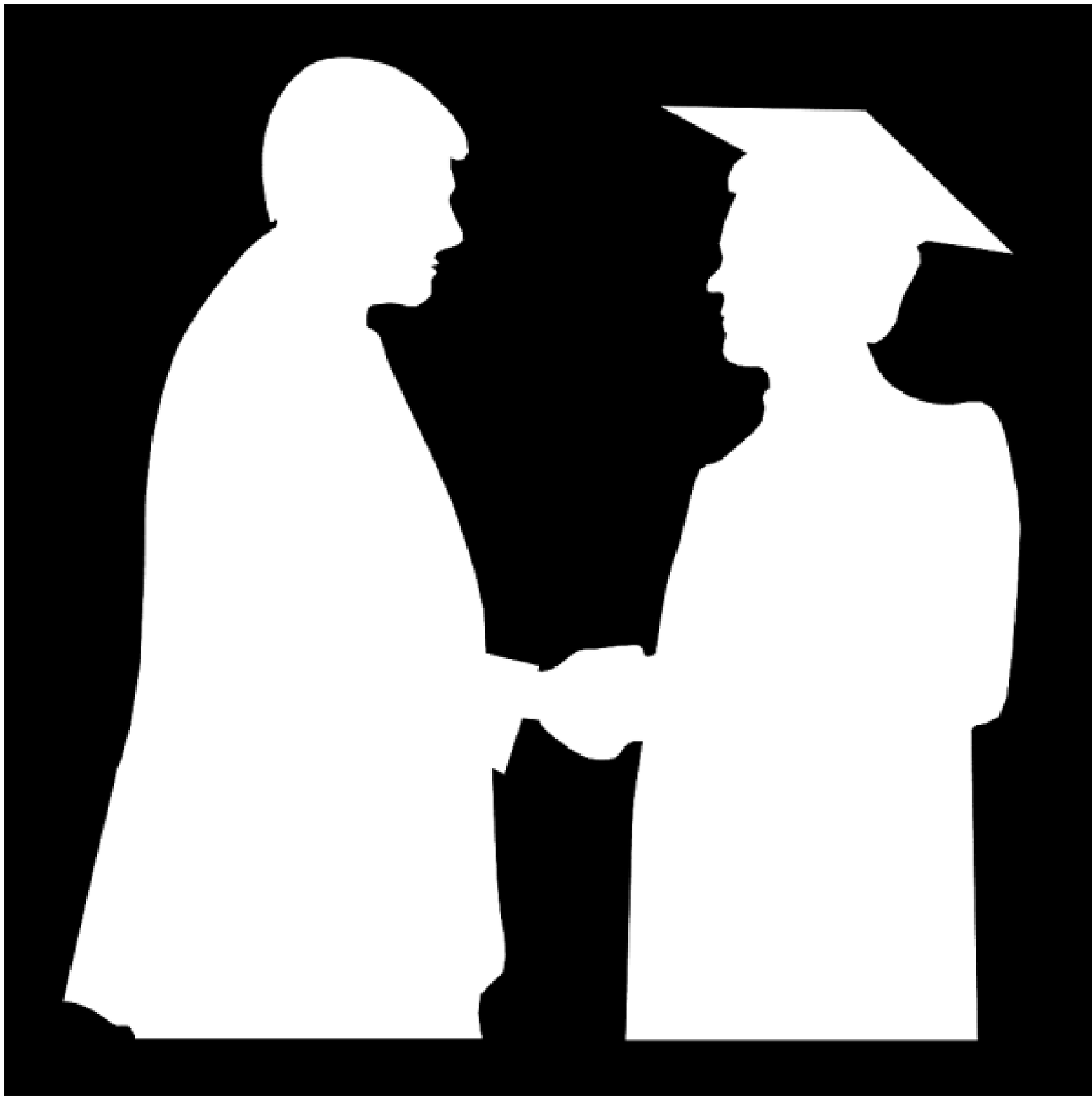}
 }
 \put(100,0){
 \includegraphics[width=0.4\textwidth]{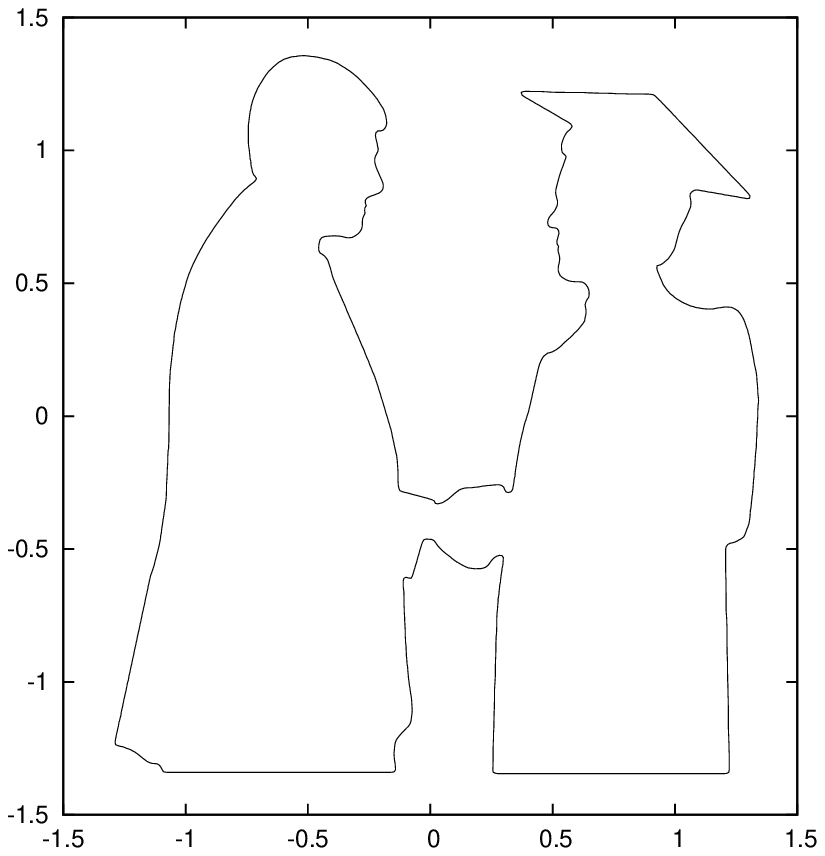}
 }
 \end{picture}
 \end{center}
\caption{The original bitmap image (left) 
and very accurate image segmentation at $t=0.025$.
We chose $N=1500$, $\varepsilon = 0$, and $F\in (-100,100)$ (right). 
The original image is taken from: {\tt http://www.zekam.uni-bremen.de/siluette.gif}. }
\label{fig:p_bmp}
\end{figure}

In this section, we will describe how to use our method for the image segmentation.
Moreover, we will compare the uniform and curvature adjusted redistribution of points.
For our computations, we did not use real images but artificial ones 
as they are better for explaining the properties of curve evolution and
points redistribution. We have selected three images for this article. 
The first one (Figure \ref{fig:yasu_bmp}) is a kanji (Chinese character) 
from Professor Mimura's name. 
You can see his name in kanji in Figure \ref{fig:mimura}. 
The shape of the character `yasu' is a relatively complex
because it contains  many gaps and sharp corners. 
The second example image (Figure \ref{fig:s_bmp}) contains 
a lot of noise and the character `S' is damaged by it. 
Finally, the third image (Figure \ref{fig:p_bmp}) contains a silhouette of
two people for which it is difficult to segment correctly the shape of their faces.

In all following examples, the initial condition is a circle with the radius $1.5$. 
Then, according to the algorithm described in the previous
section, we let the curve evolve until the stopping condition is 
reached. It is important to pay attention to proper choice of 
segmentation parameters.

The image segmentation by this method has three important parameters. 
The parameter $N$ denotes the number of points  which the curve consists of.
It is clear that the more points we use the better final shape we get. 
On the other hand, higher number of grid points makes the computation much slower.
The parameter $\varepsilon$ (see equation (\ref{eq:phi})) controls the curvature adjusted redistribution. 
For $\varepsilon=0$ we obtain the asymptotically uniform redistribution. 
For $\varepsilon\in(0,1)$ we obtain a curvature adjusted redistribution. 
For numerical computations, we recommend to choose the 
value between $0.1$ and $0.2$. The last but very important function is 
the external force term $F$ depending on given $600\times 600$ pixels bitmap picture 
and influences the final shape significantly. 

In the following computations, 
the target figure is given by a digital gray scale bitmap image 
represented by integer values between $0$ and $255$ on $600\times 600$ pixels. 
The values $0$ and $255$ correspond to black and white colors, respectively, 
whereas the values between $0$ and $255$ correspond to gray colors. 
In our three examples, the target shapes are given in white
color with black background.

Given a figure, we can construct its
image intensity function $I(\vecx)\in \{0, \ldots, 255\}\subset\Z$ 
defined in the computational domain
$\Omega:=(-1.5,\,1.5)\times (-1.5,\,1.5)$. 
We remark that $I(\vecx)$ is 
piecewise constant in each pixel.

\begin{figure}[ht]
 \begin{center}
 \subfigure[$F\in(-20,30)$]{
 \label{fig:yasu_f20}
 \includegraphics[width=0.22\textwidth]{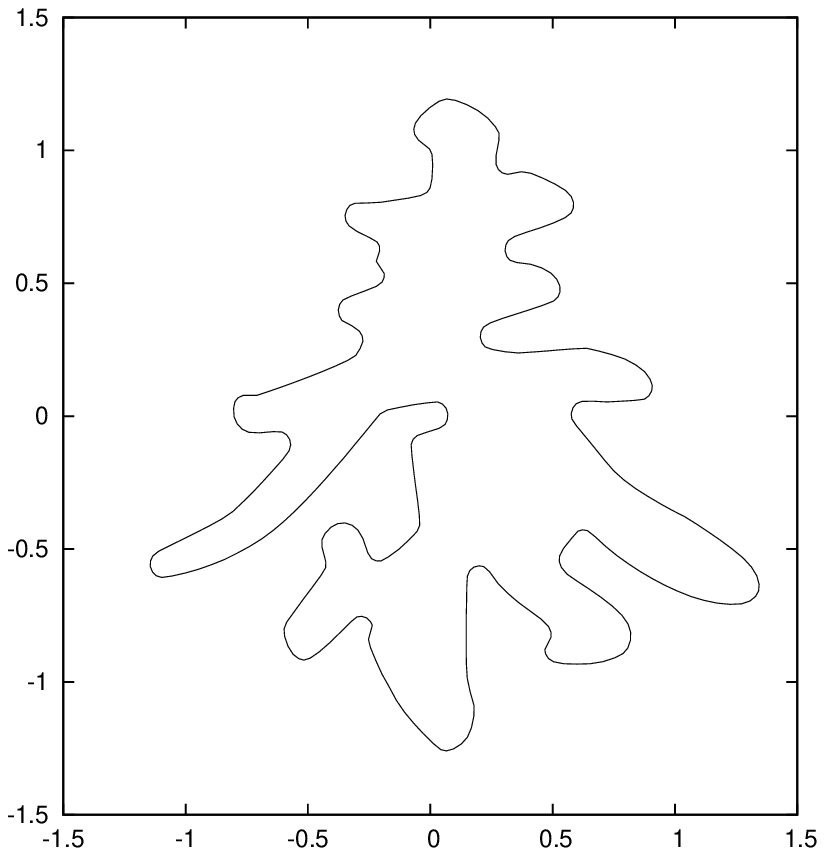}
 }
\subfigure[$F\in(-30,30)$]{
 \label{fig:yasu_f30}
 \includegraphics[width=0.22\textwidth]{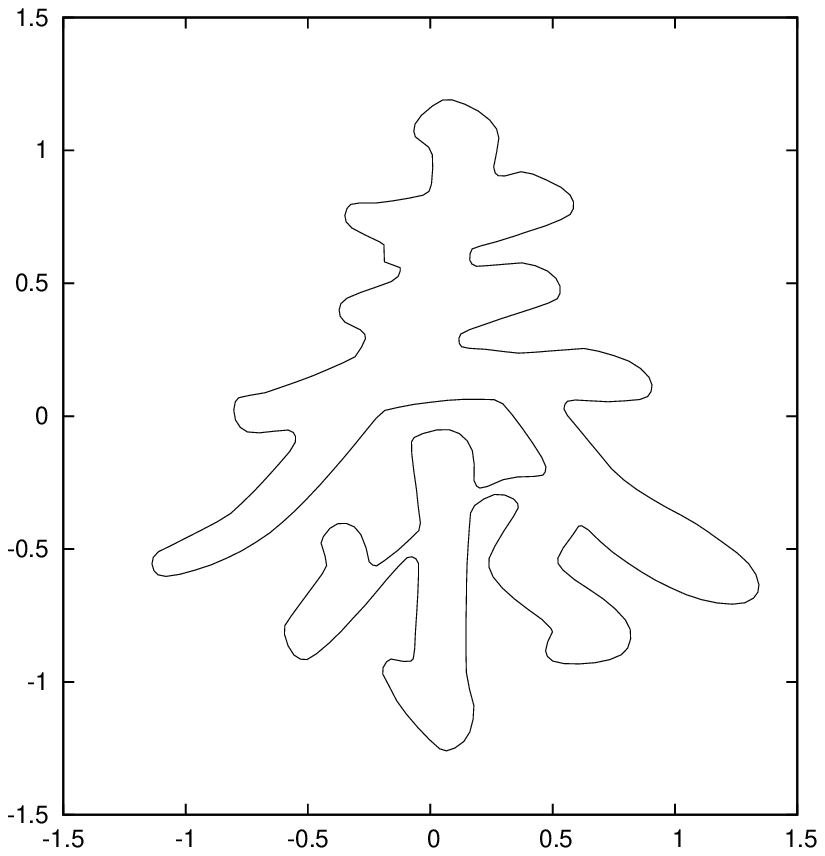}
 }
\subfigure[$F\in(-40,40)$]{
 \label{fig:yasu_f40}
 \includegraphics[width=0.22\textwidth]{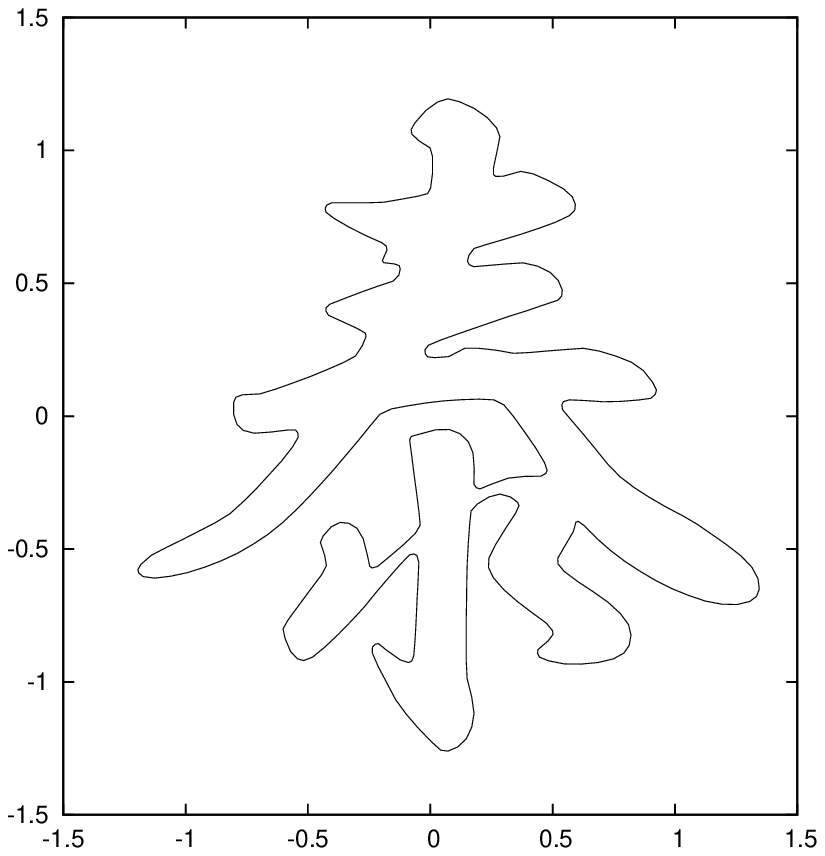}
 }
\subfigure[$F\in(-100,100)$]{
 \label{fig:yasu_f100}
 \includegraphics[width=0.22\textwidth]{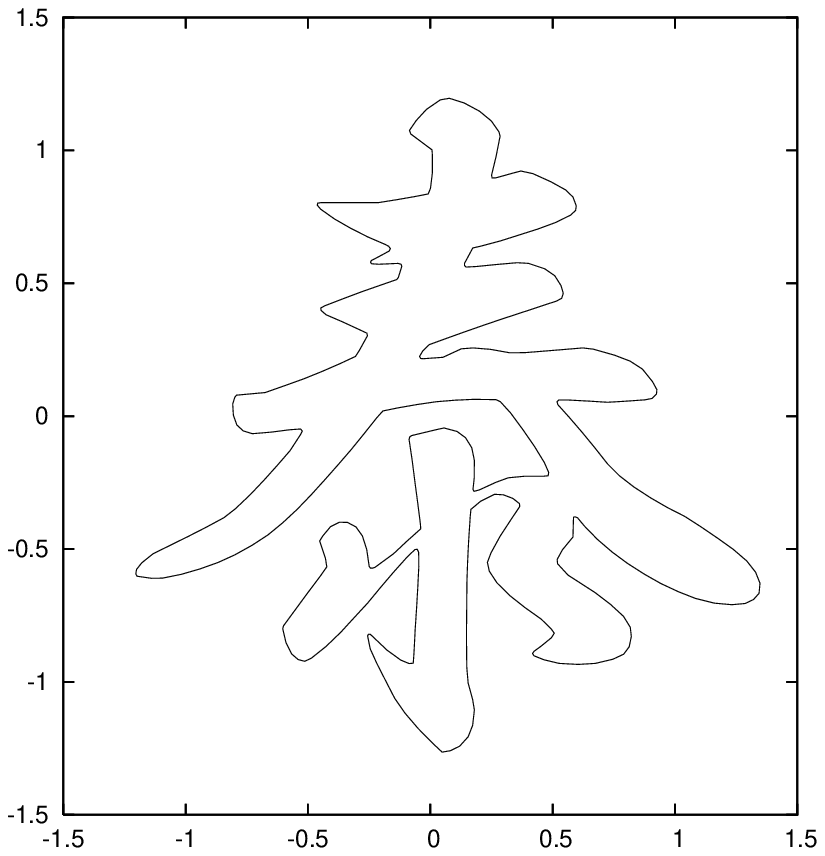}
 }
 \end{center}
\caption{The comparison of several choices of the external force
with $N=350$ and $\varepsilon = 0.2$.}
\label{fig:force}
\end{figure}

We define the forcing term $F(\vecx)$ as follows:
\begin{equation}
\label{eq:color}
F(\vecx)=F_{max}-(F_{max}-F_{min})\frac{I(\vecx)}{255}
\quad(\vecx\in\Omega),
\end{equation}
where $F_{max}>0$ corresponds to purely black color (background) and $F_{min}<0$
corresponds to purely white color (the object to be segmented). Maximal and
minimal values determine the final shape because in general $1/F$ is equivalent to
the minimal radius the curve can attain. The choice of small values
of $F$ causes the final shape to be rounded or the curve can not pass through
narrow gaps. Various choices of $F_{max}$ and $F_{min}$ are shown in Figure~\ref{fig:force}. 
Figure~\ref{fig:yasu_f20} illustrates a
situation when the force is too small and the curve can not evolve through
the obstacle.
Increasing the value of $F$ 
(Figures~\ref{fig:yasu_f30}, \ref{fig:yasu_f40}, \ref{fig:yasu_f100}), 
the shape becomes more and more sharp. 

\begin{figure}[ht]
\begin{center}
 \subfigure[$t=0.005$]{
 \includegraphics[width=0.17\textwidth]{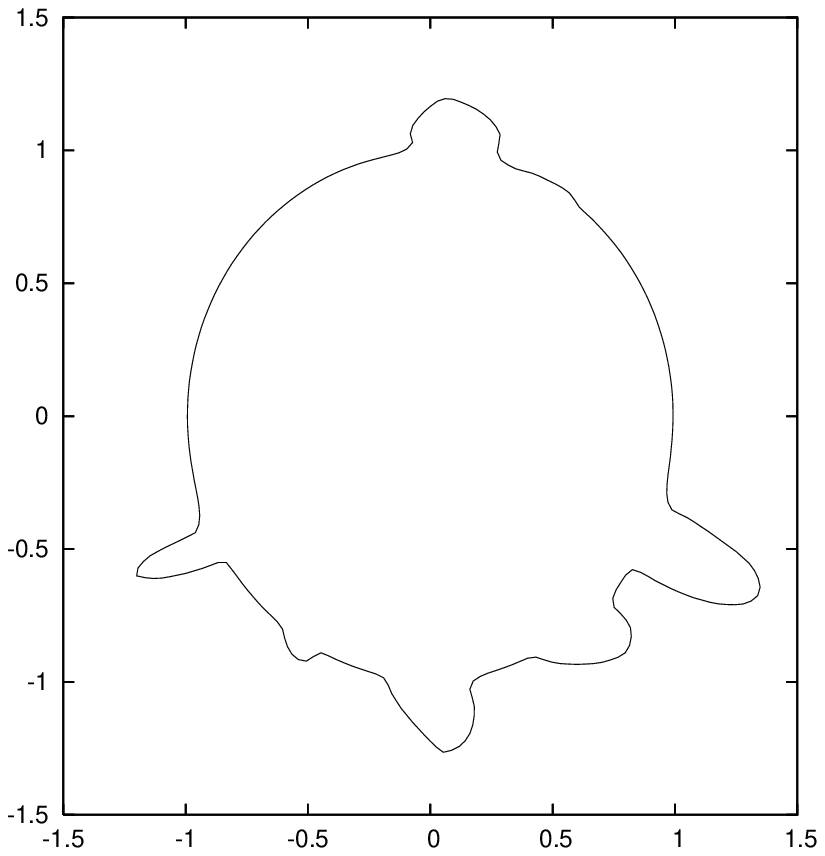}
 \label{fig:yasu350e0_1}
 }
 \subfigure[$t=0.010$]{
 \label{fig:yasu350e0_2}
 \includegraphics[width=0.17\textwidth]{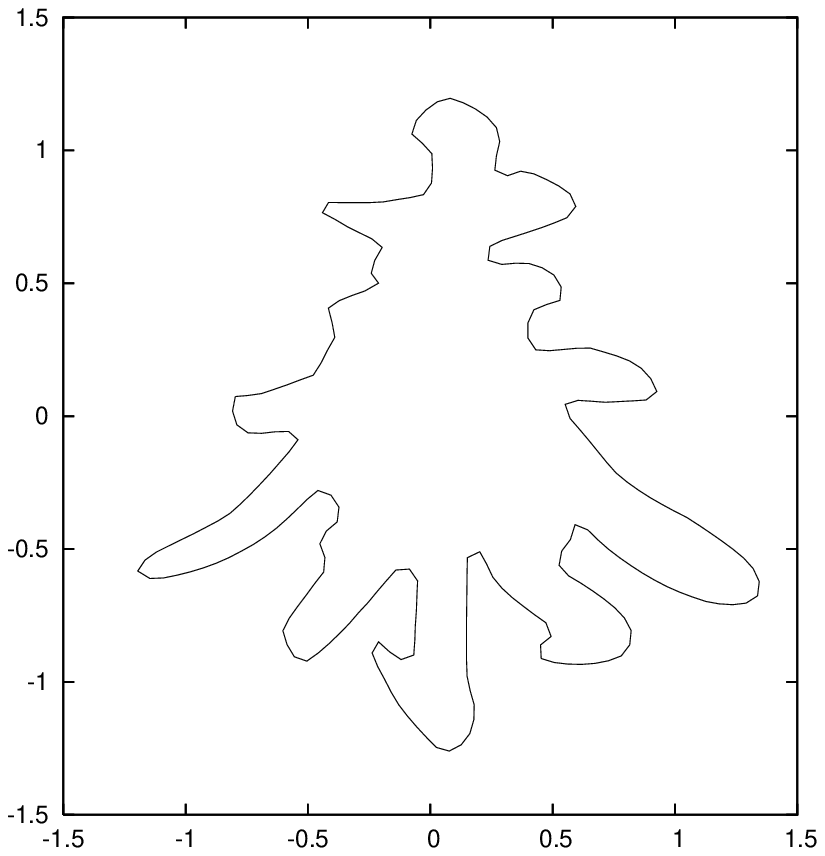}
 }
 \subfigure[$t=0.015$]{
 \label{fig:yasu350e0_3}
 \includegraphics[width=0.17\textwidth]{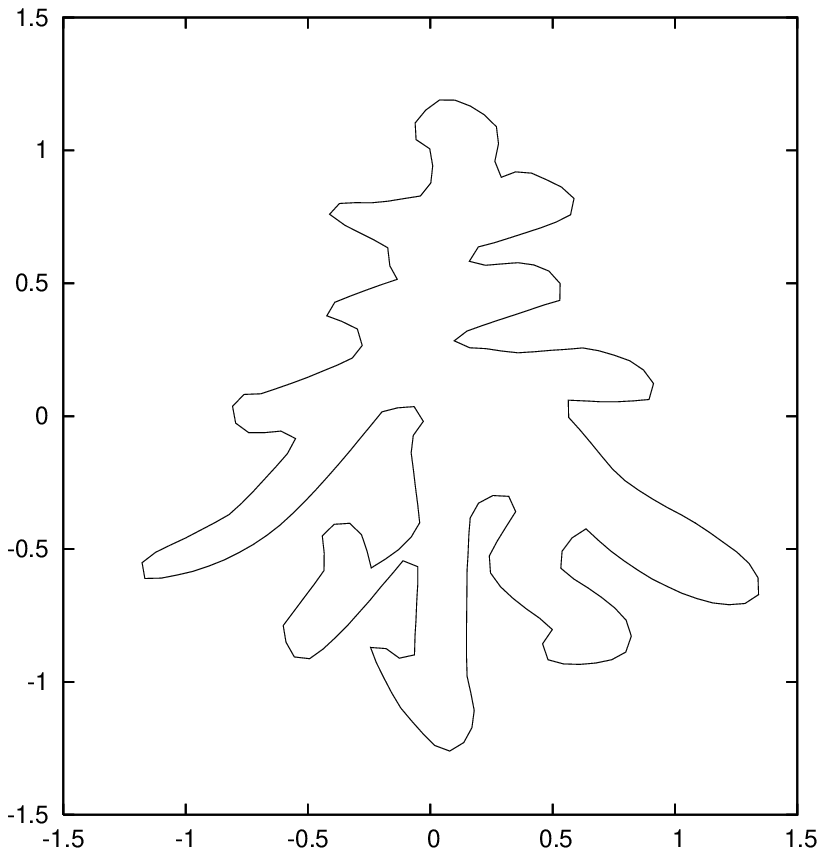}
 }
 \subfigure[$t=0.020$]{
 \label{fig:yasu350e0_4}
 \includegraphics[width=0.17\textwidth]{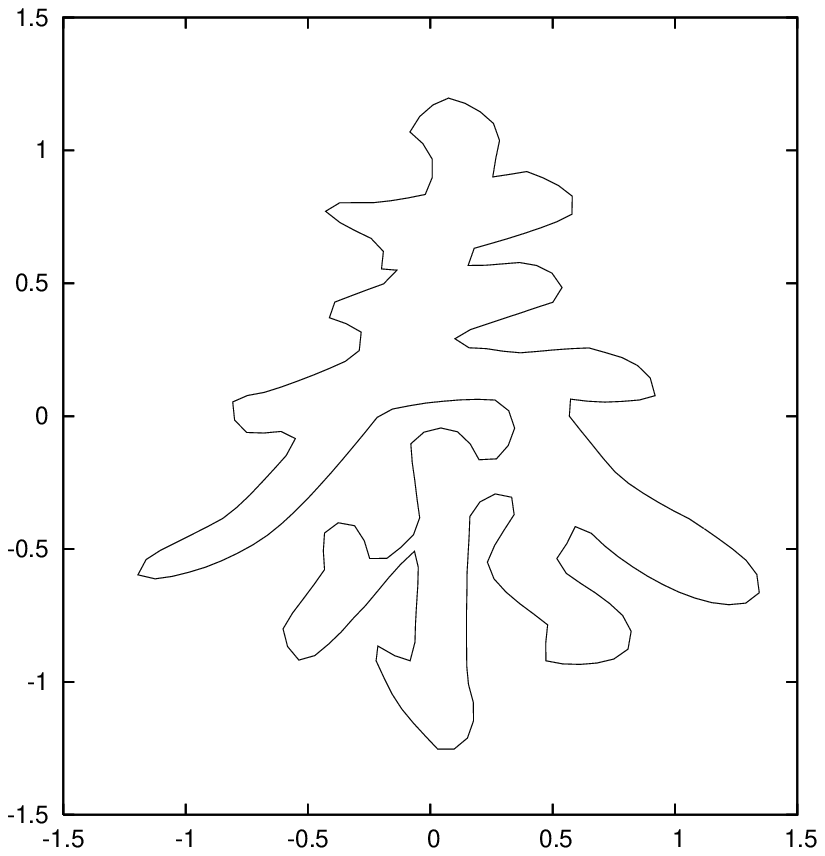}
 }
 \subfigure[$t=0.025$]{
 \label{fig:yasu350e0_5}
 \includegraphics[width=0.17\textwidth]{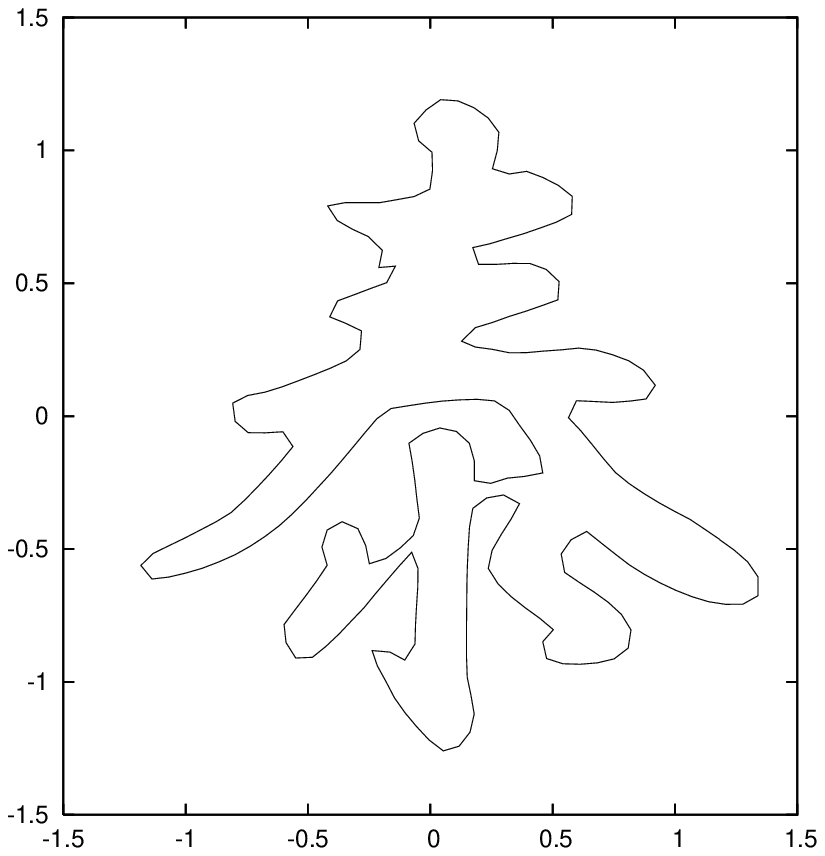}
 }
\end{center}
\caption{$N=250$, $\varepsilon = 0$, $F\in[-100,100]$.}
\label{fig:yasu_eps0}
\end{figure}

\begin{figure}[ht]
 \begin{center}
 \subfigure[$t=0.005$]{
 \label{fig:yasu350e02_1}
 \includegraphics[width=0.17\textwidth]{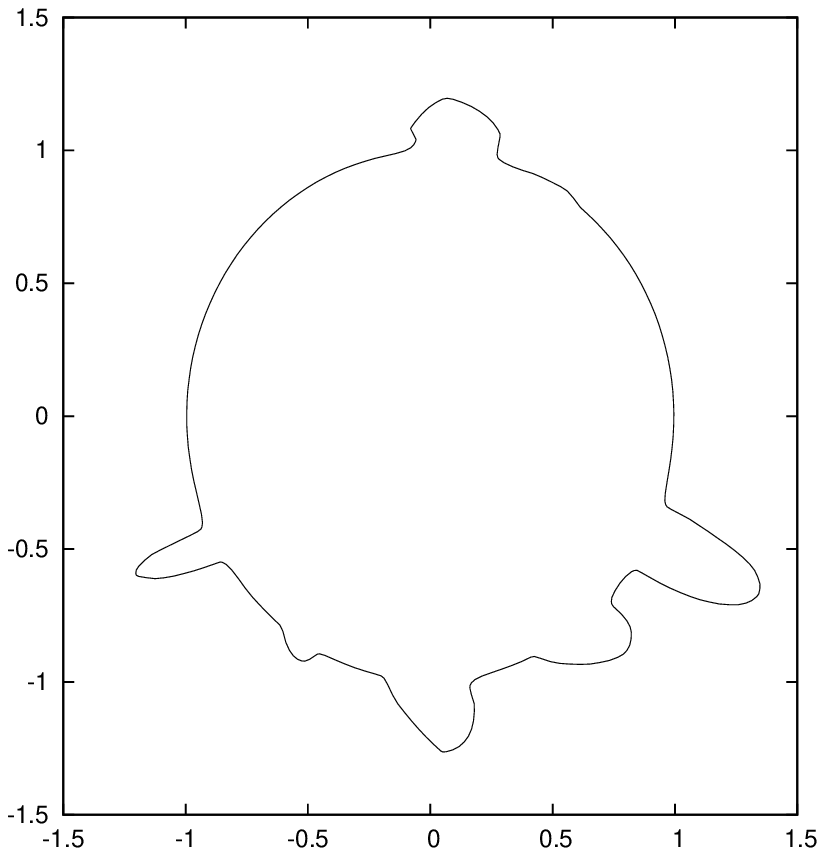}
 }
 \subfigure[$t=0.010$]{
 \label{fig:yasu350e02_2}
 \includegraphics[width=0.17\textwidth]{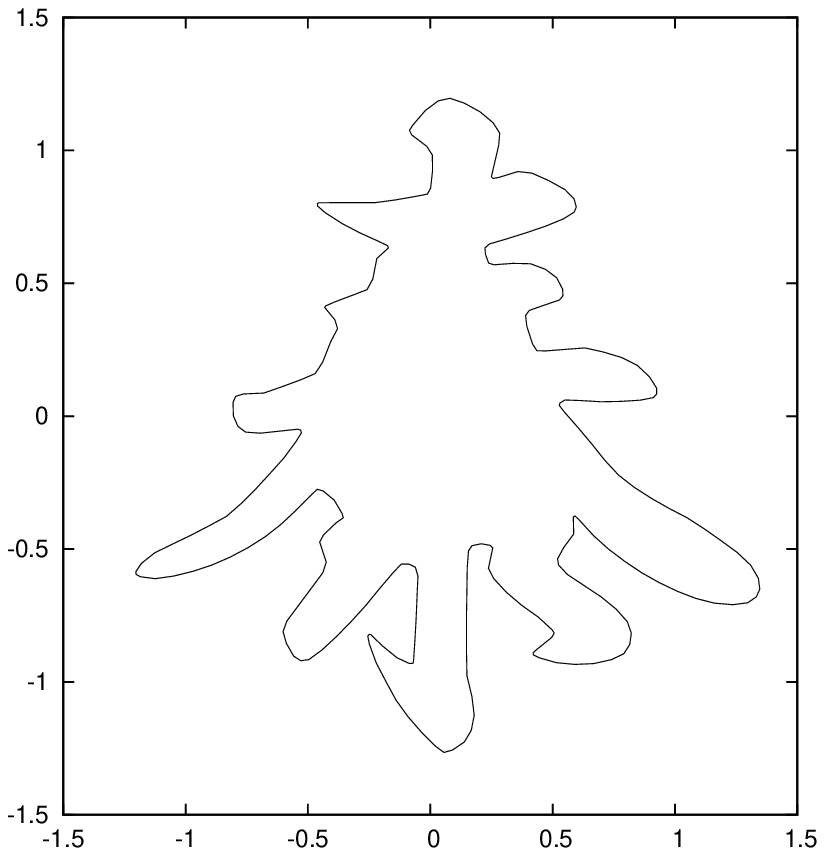}
 }
 \subfigure[$t=0.015$]{
 \label{fig:yasu350e02_3}
 \includegraphics[width=0.17\textwidth]{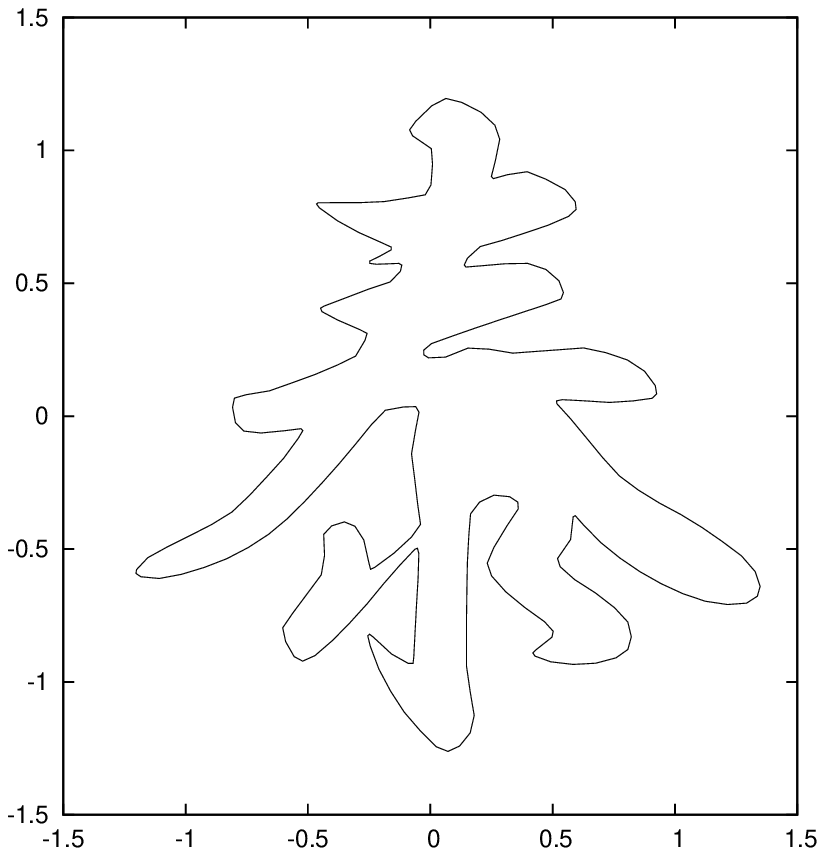}
 }
 \subfigure[$t=0.020$]{
 \label{fig:yasu350e02_4}
 \includegraphics[width=0.17\textwidth]{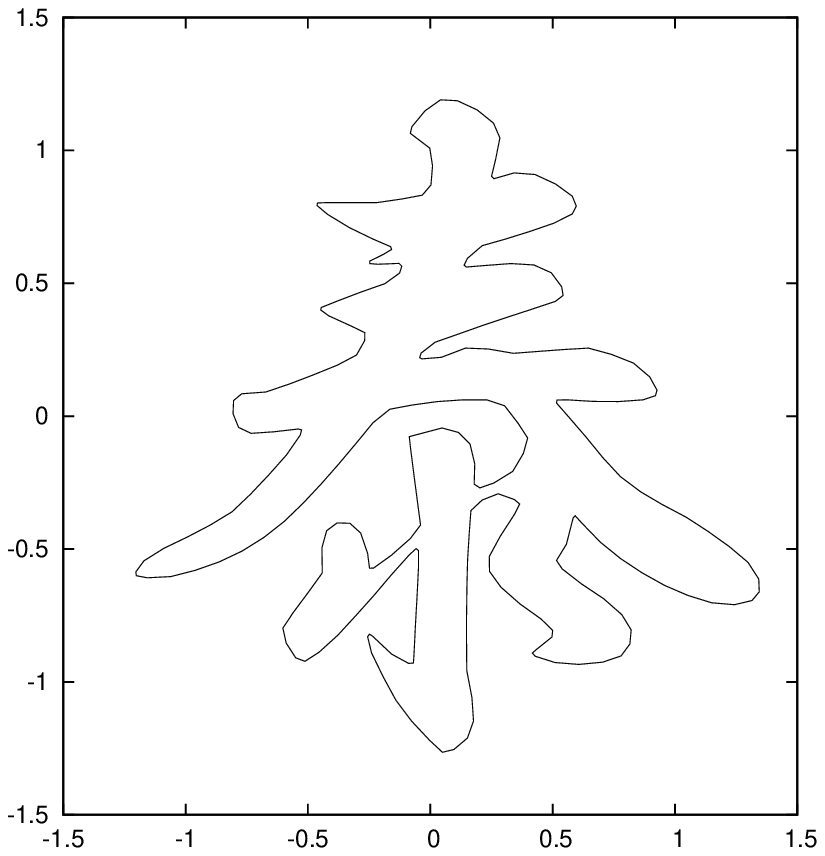}
 }
 \subfigure[$t=0.025$]{
 \label{fig:yasu350e02_5}
 \includegraphics[width=0.17\textwidth]{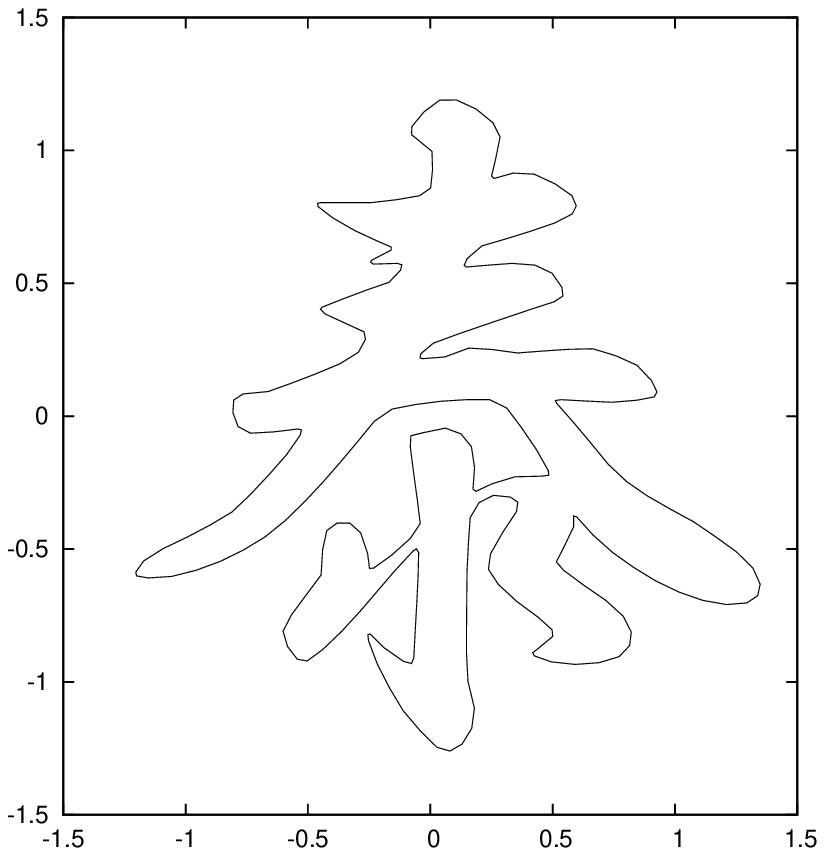}
 }
 \end{center}
\caption{$N=250$, $\varepsilon = 0.2$, $F\in[-100,100]$.}
\label{fig:yasu_eps02}
\end{figure}

Now we will describe the results of a uniform ($\varepsilon=0$) and curvature adjusted ($\varepsilon>0$) redistribution method.
Figures~\ref{fig:yasu_eps0} and \ref{fig:yasu_eps02} show the
evolution of segmentation curves in time. 
In both cases, the number of points $N=250$, 
the external forcing term $F$ is calculated by (\ref{eq:color}) where 
$F_{min}=-100$ and $F_{max}=100$ but the $\varepsilon$ is different. One can see that for $\varepsilon=0$ (Figure \ref{fig:yasu_eps0}) very sharp corners are smoothed or contain too few points 
while for $\varepsilon=0.2$ (Figure \ref{fig:yasu_eps02})
even the spike under the first horizontal line is kept. Relatively high
value of the forcing term makes the evolution fast. Both computations were 
performed for the same time interval $t\in(0, 0.025)$. 
Figure~\ref{fig:yasu_bmp}(right) shows the result with very high number points $N=1500$
and $\varepsilon=0$.

\begin{figure}[ht]
 \begin{center}
 \subfigure[$N=250$]{
 \label{fig:yasu_comp250}
 \includegraphics[width=0.4\textwidth]{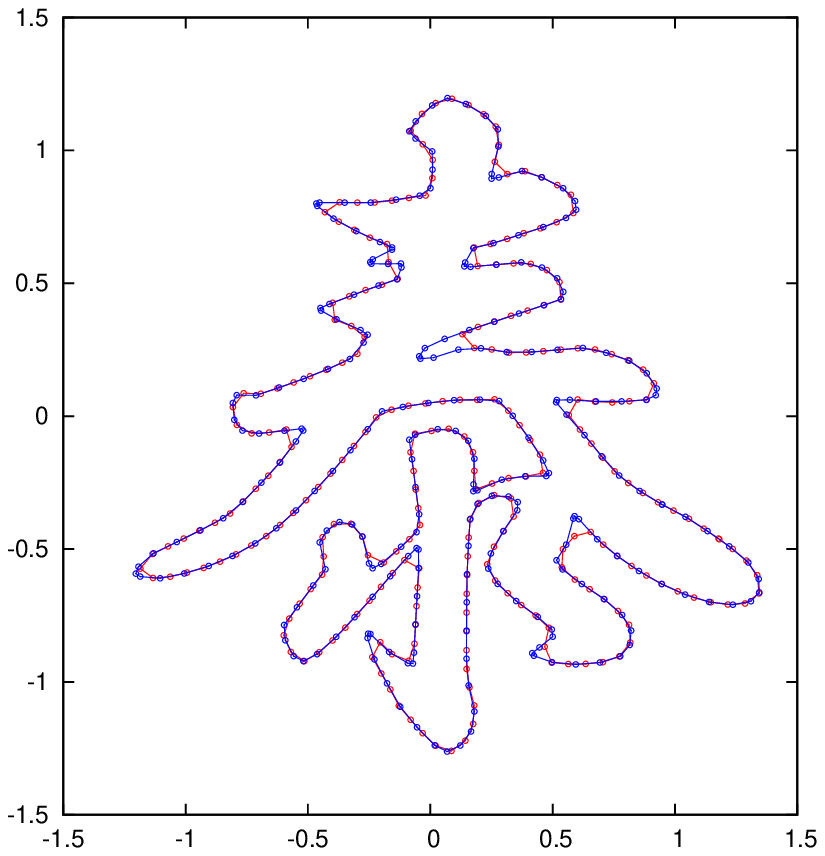}
 }
 \subfigure[$N=350$]{
 \label{fig:yasu_comp350}
 \includegraphics[width=0.4\textwidth]{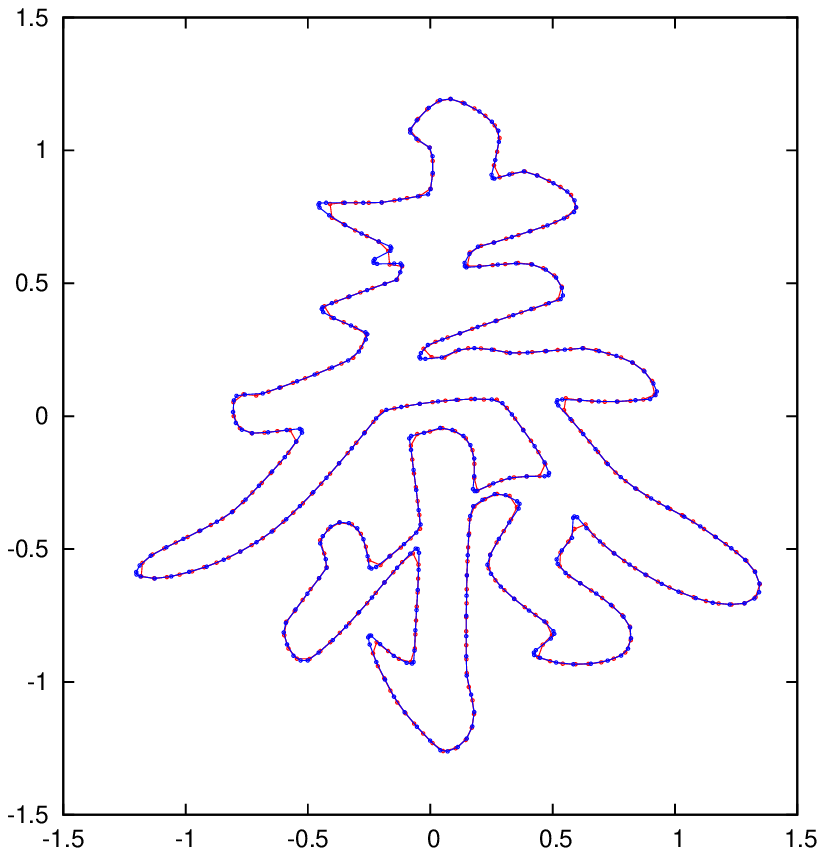}
 }
 \end{center}
\caption{A comparison between image segmentations with $\varepsilon=0$ and $\varepsilon=0.2$.
Here $F\in(-100,100)$ and $t=0.025$.}
\label{fig:yasu_comp}
\end{figure}

Figure~\ref{fig:yasu_comp} illustrates a comparison 
of the last step of the evolution for the uniform and
curvature adjusted redistribution for two different choices of $N$.
The curve in Figure \ref{fig:yasu_comp250} contains $250$ grid points.
One can clearly see that for $\varepsilon=0$ the narrow gap under
the second horizontal line is not segmented well and also 
high curvature parts are rounded. 
For $N=350$ (Figure~\ref{fig:yasu_comp350}), the segmentation is quite
similar for both cases but still for $\varepsilon=0$, 
sharp parts are not as good as for the case $\varepsilon=0.2$. 

We tested the redistribution method also for images containing noise 
or some additional artifacts. Figure~\ref{fig:s_bmp} illustrates this case. 
According to (\ref{eq:color}), gray colors do not generate so strong force.
Hence the curve can easily pass through. If the noise generates too large 
external force (almost white color) then it prevents the curve to move further. 
But on the 
other hand, very high curvature will appear and causes the curve to overcome the noise. 
When we work with noisy images, it is very important
to choose suitable interval for $F$. If too wide interval is chosen then
even a small noise can stop the computation. 
On the other hand, too narrow interval will cause bad shape segmentation. 
In our computations for noisy images, we chose $F\in (-30,35)$. 

\begin{figure}[ht]
 \begin{center}
 \subfigure[$N=80$]{
 \label{fig:s_comp080}
 \includegraphics[width=0.4\textwidth]{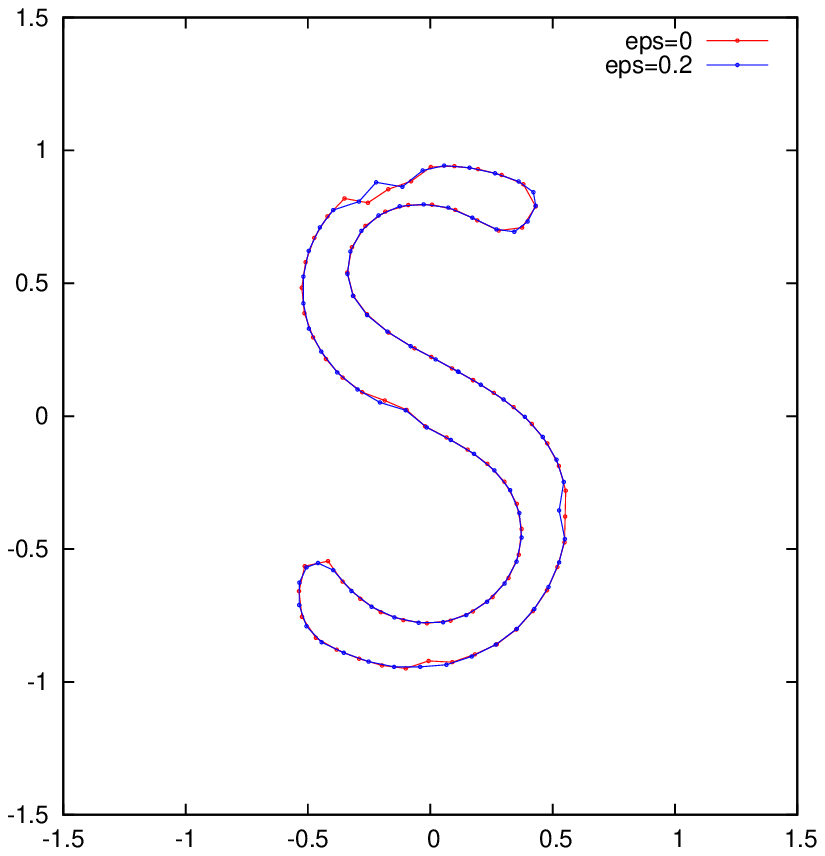}
 }
 \subfigure[$N=150$]{
 \label{fig:s_comp150}
 \includegraphics[width=0.4\textwidth]{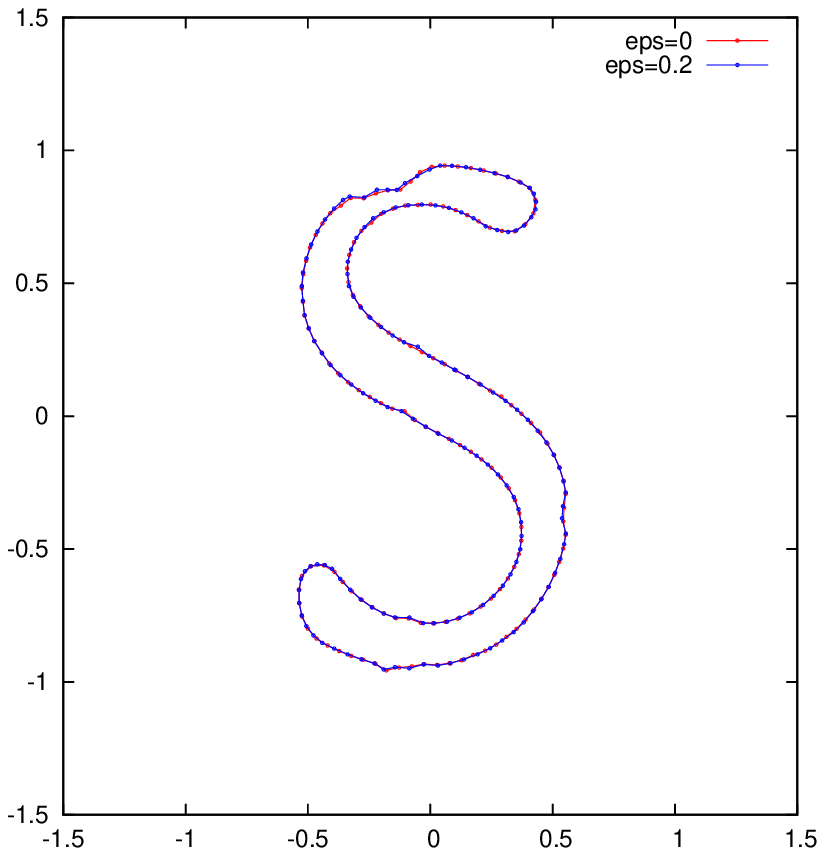}
 }
 \end{center}
\caption{A comparison between image segmentations
with $\varepsilon=0$ and $\varepsilon=0.2$. Here $F\in(-30,35)$ and $t=0.1$.}
\label{fig:s_comp}
\end{figure}

Since the shape in Figure~\ref{fig:s_bmp} does not have many details
and the length is not so big, we do not need to use high number of points
$N$. The computation in Figure~\ref{fig:s_bmp}(right) was done for $N=800$ 
only for a reference value. The comparison has been performed for $N=80$ 
(Figure~\ref{fig:s_comp080} and $N=150$ (Figure~\ref{fig:s_comp150}). 

\begin{figure}[ht]
 \begin{center}
 \subfigure[$N=150$]{
 \label{fig:p_comp150}
 \includegraphics[width=0.4\textwidth]{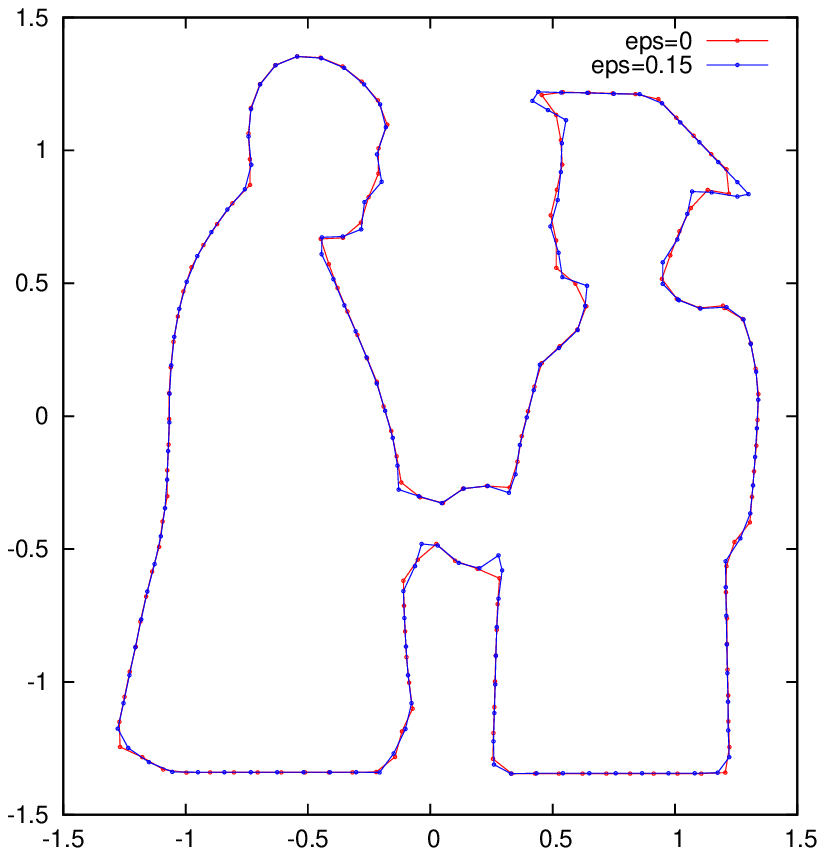}
 }
 \subfigure[$N=350$]{
 \label{fig:p_comp350}
 \includegraphics[width=0.4\textwidth]{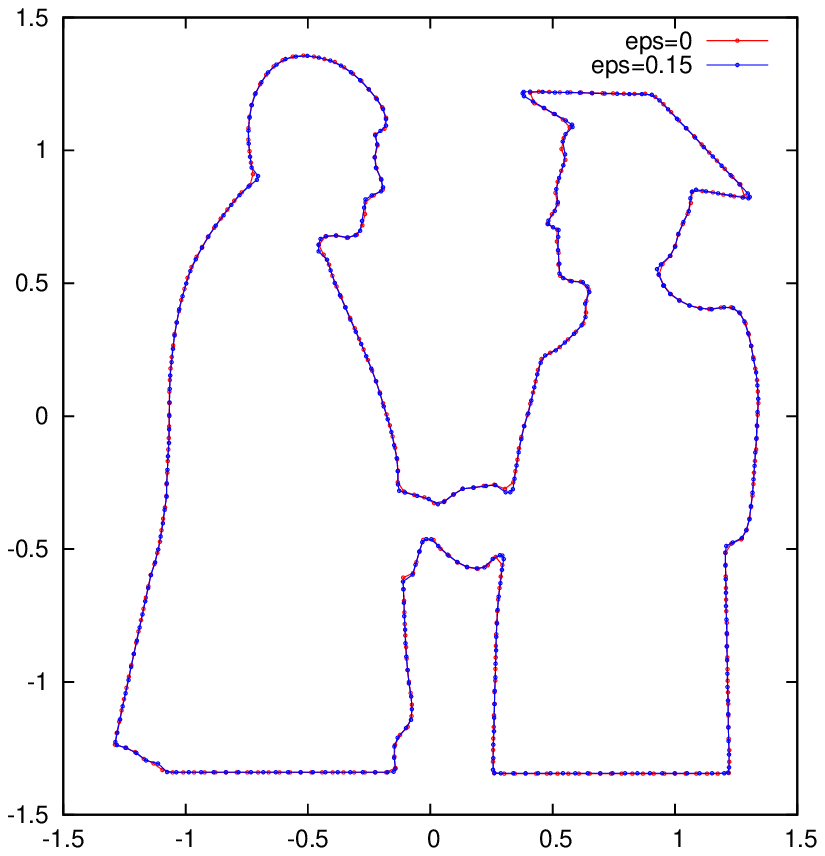}
 }
 \end{center}
\caption{A comparison between image segmentations with $\varepsilon=0$ and $\varepsilon=0.15$. 
Here $F\in(-100,100)$ and $t=0.05$.}
\label{fig:p_comp}
\end{figure}

The last example depicts silhouette of two persons (Figure \ref{fig:p_bmp}). 
There are parts having high curvature and the faces of persons have small details. 
We need a large forcing term to find the shape correctly. 
Since there is no noise, we can do it. 
Even the simulation for a quite low $N=150$, 
Figure \ref{fig:p_comp150} shows that we still get a reasonably
good approximation with $\varepsilon=0.15$. Computation with 
$\varepsilon=0$ did not find details near hands and also hat is not
as sharp as it should be. By increasing $N$ to 350, we were able to achieve a good 
approximation with both values of $\varepsilon$.
But again, sharp edges
in the hat area are better segmented with $\varepsilon=0.2$ (Figure \ref{fig:p_comp350}). 

All computations were performed on a standard personal computer with
Intel processor at 2.4 GHz. The time of computation was never higher
than about 10 seconds. The computation was faster for $\varepsilon=0$ 
than for $\varepsilon>0$. For $N$ about 200, the CPU time is usually
less than 1 second.

\end{document}